%% file: article.tex
\documentclass[onefignum,onetabnum]{siamart251216}

\input{shared}

\ifpdf
\hypersetup{
  pdftitle={Local optimization of weak distance between compact surfaces on special Euclidean group},
  pdfauthor={K. Koga}
}
\fi

\DeclareFontFamily{U}{mathx}{}
\DeclareFontShape{U}{mathx}{m}{n}{<-> mathx10}{}
\DeclareSymbolFont{mathx}{U}{mathx}{m}{n}
\DeclareMathAccent{\widecheck}{0}{mathx}{"71}

\begin{document}

\maketitle

\begin{abstract}
We consider local optimization of a weak distance between two compact surfaces embedded in the three-dimensional Euclidean space on its special Euclidean group. Identifying those objects with the associated surface measures, their discrepancy is quantified in terms of the inhomogeneous Sobolev norm of negative order via the Plancherel theorem. Then, applying an isometry to one surface corresponds to pushforwarding its surface measure and the distance can be regarded as a function on the Lie group. For appropriate exponents of the Sobolev norm, the second power of the function acquires sufficient differentiability that allows to search for its local minima in a derivative-based framework, and the gradient of the objective function has a favorable structure for efficient implementations using the nonuniform fast Fourier transform. In numerical experiments, we observe convergence of the SR1 trust-region method applied to a few root-finding problems and discuss its connection to a more geometric quantity.
\end{abstract}




\section{Introduction}\label{sec1}
Geometric data have been rapidly accumulating in the field of biology, chemistry, and others where observations using electron microscopy play a crucial role. This trend is emphasized by, for example, the emergence of public databases such as the Electron Microscopy Data Bank (EMDB) for EM density maps obtained from cryo-electron microscopy (Cryo-EM) in protein structure analysis, for which a number of solutions to alignment problems have been proposed \cite{SiYa2024, RiZhCoChDu2025, RiWoPoCoDu2023, PeGoHuMeCoCrMoFe2021, HaTeChChKi2021}. On the other hand, the recent advancements in generating surface models of embryonic tissues from images acquired via the confocal fluorescence microscopy \cite{IcSiDeTu2025, IcDeMcDuTu2023} can be a strong motivation toward developing analogues of tools for EM density maps, whereas its progress is behind the case of such function data. \par

In general, alignment problems are classified into two types: global and partial. Loosely speaking, global alignment determines whether two geometric data are the same as a whole while partial alignment judges if one is part of the other. In either case, it is often insufficient to define a naive distance between given objects, because a computationally tractable quantity may lead one to a wrong conclusion that two data identical up to isometry are different. To address this issue for surfaces, several authors have suggested intrinsic methods that essentially utilize information on the objects only. For instance, Bronstein, Bronstein, and Kimmel \cite{BrBrKi2006} develops practical algorithms for computing the Gromov-Hausdorff distance \cite{BrHa1999} between triangulated surfaces, which measures discrepancies between metric spaces uniformly, and its asymmetric variant called the partial embedding distance. In another direction, Lipman and Daubechies \cite{LiDa2011} proposes polynomial-time schemes based on the conformal Wasserstein distance, which reduce comparisons between disc-like surfaces to the Kantorovich problem \cite{Villani2003} on the unit disk via conformal mappings, and the strategy is applied to evaluating the geometric similarity of anatomical surfaces \cite{BoLiClPuPaFuJeDa2011}. However, the former relies on exhaustive search for optimal embeddings between surfaces and hence leads to high-dimensional optimization problems, while the latter is topology-dependent and difficult to extend to general cases.\par

More recently, alignment of protein data in the format of EM density maps has served as a strong driving force towards developing isometry-invariant distances, and several extrinsic algorithms based on optimal transport have been suggested in the research domain. Those include, but are not limited to, bayesian optimization of the wavelet earth mover's distance \cite{SiYa2024} and an application of the (unbalanced) Gromov-Wasserstein theory \cite{RiZhCoChDu2025}. While such extrinsic concepts can be generalized to comparisons between surfaces \cite{GoNiLiCrSaCuZoHaCa2023}, however, practical usage of the existing methods is to generate rough approximations of the desired minimizers, and typically alignment must be completed by a local optimization technique (called refinement in the literature) such as the \texttt{fitmap} command in ChimeraX \cite{PeGoHuMeCoCrMoFe2021} that maximizes correlations between volume data via the steepest gradient ascent or the local search in AlignOT \cite{RiWoPoCoDu2023} that employs the stochastic gradient descent for minimizing an entropic regularization of the quadratic Wasserstein distance. \par

As another extrinsic formulation, Koga \cite{Koga2024} introduces a weak distance between two compact surfaces embedded in the three-dimensional (3D) Euclidean space $\mathbb{R}^3$. Mathematically, it identifies given surfaces with the associated surface measures and quantifies their discrepancy as elements of the inhomogeneous Sobolev space $H_s(\mathbb{R}^3)$ for negative $s$, where the norm $\|\cdot \|_{H_s(\mathbb{R}^3)}$ can be directly evaluated from their Fourier transforms via the Plancherel theorem. Here, compactness guarantees smoothness of the Fourier transforms, and the Sobolev norm can be well approximated by the trapezoidal rule on a bounded domain. Furthermore, surface measures are also discretized with fast-converging numerical quadratures on simplices, and subsequent computations of the Fourier transforms are efficiently performed by the combination of the nonuniform fast Fourier transform (NUFFT) \cite{BaMaKl2019, ShWrAnBlBa2021} and parallelization on a single Graphics Processing Unit (GPU). As a result, the suggested algorithm is fast and accurate, and it is iteration-free as opposed to those for the Hausdorff and Wasserstein distances. However, the quantity depends on how we embed two surfaces in the ambient space, and it is also nontrivial that such a weak concept retains geometric information that is useful in applications, while there are evidences that the truncated Fourier transforms of singular measures may provide a solid basis in geometry processing \cite{JiWaHuMaNi2019,IcSiDeTu2025}.  \par

Having those questions in mind, this paper develops a novel algorithm for performing local optimization of the weak distance between compact surfaces \cite{Koga2024} on the special Euclidean group $SE(3)$. Namely, by applying an isometry to one surface, which corresponds to pushforwarding its surface measure, we regard the distance as a function on the space of isometries. For appropriate exponents of the Sobolev norm, the second power of the function acquires sufficient differentiability that allows to search for its local minima in a derivative-based framework. Moreover, it turns out that essential part of computing the gradient of the objective function is reduced to evaluations of three additional Fourier-type transforms, each of which can be efficiently approximated by NUFFT. In numerical experiments, we apply the trust-region method with the symmetric rank-1 (SR1) Hessian updates and confirm its convergence in a few root-finding problems for triangulated surfaces. We experimentally observe that the low dimensionality of the optimization problem enables the suggested algorithm to accurately return the global minimizer provided the initial guess is appropriate, and also verify that the maximum value of the Euclidean distance between the nodes constituting triangulated surfaces closely follows behavior of the objective function, which implies a tight link between the weak distance and more geometric quantities.\par

The rest of this paper is organized as follows. Section \ref{formulate} reviews the distance between compact surfaces in \cite{Koga2024} and formulates its local optimization problem on the special Euclidian group. Section \ref{method} explains numerical methods and some techniques in their efficient implementations. Section \ref{results} presents numerical results and evaluates the implemented codes. Concluding remarks are given in Section \ref{conclusion}.


\section{Formulations}\label{formulate}
This section formulates the optimization problem of the weak distance between compact surfaces in terms of the inhomogeneous Sobolev norm \cite{Koga2024} on the special Euclidean group $SE(3)$. To this end, we first review basic concepts in Euclidean harmonic analysis \cite{Grafakos2014-1,Grafakos2014-2}, which is followed by introducing the objective function and deriving its gradient via the elementary Lie theory.
%
%
%
%
%
%
%
\subsection{Fourier transform on $\mathbb{R}^3$}
Given a Borel-measurable function $f$ on $\mathbb{R}^3$, its Fourier transform $\widehat{f}$ is defined as the integrals
\begin{equation}
\label{eq:ft_forward}
\widehat{f}(\xi) = \int_{\mathbb{R}^3} f(x) e^{-2\pi i \xi \cdot x} dx,\quad \xi \in \mathbb{R}^3,
\end{equation}
where $dx$ denotes the Lebesgue measure and we call the variable $\xi$ the frequency. For a wide class of functions, the inversion formula
\begin{equation}
\label{eq:ft_inverse}
\widecheck{f}(x) = \int_{\mathbb{R}^3} f(\xi) e^{2\pi i \xi \cdot x} d\xi,\quad x \in \mathbb{R}^3,
\end{equation}
recovers the original $f$ from (\ref{eq:ft_forward}) in a suitable sense. As an extension of (\ref{eq:ft_forward}), the Fourier transform of a Borel measure $\mu$ is defined as
\begin{equation}
\label{eq:ft_measure}
\widehat{\mu}(\xi) = \int_{\mathbb{R}^3}  e^{-2\pi i \xi \cdot x} d\mu(x),\quad \xi \in \mathbb{R}^3,
\end{equation}
whereas the inversion (\ref{eq:ft_inverse}) no longer makes sense in general. For instance, the Dirac measure $\delta_0$ at the origin is Borel and its Fourier transform is
$\widehat{\delta_0}\equiv 1$, which only belongs to the class of bounded functions $L^{\infty}(\mathbb{R}^3)$. Nevertheless, the Fourier transform is further extended to elements of the space of tempered distributions $\mathcal{S}'(\mathbb{R}^3)$, which is the space of continuous linear functionals on the Schwartz space $\mathcal{S}(\mathbb{R}^3)$ and includes all finite Borel measures, and the inverse operation to (\ref{eq:ft_measure}) is defined more formally in the sense of tempered distributions.\par

One of the most important properties of the Fourier transform is that it is an isometry on the space of square-integrable functions $L^2(\mathbb{R}^3)$. Namely, if $f \in L^2(\mathbb{R}^3)$, the norm $\|f\|_{L^2(\mathbb{R}^3)}$ is computed from $\widehat{f}$ via the Plancherel theorem
\begin{equation}
\label{eq:plancherel}
\|f\|^2_{L^2(\mathbb{R}^3)}=\int_{\mathbb{R}^3} |f(x)|^2 dx=\int_{\mathbb{R}^3} |\widehat{f}(\xi)|^2 d\xi.
\end{equation}
To define a similar norm for a class of Borel measures via the the Fourier transform, one can consider the (inhomogeneous) Sobolev space, denoted by $W_s^p(\mathbb{R}^3)\, (1<p<\infty,\,s\in \mathbb{R})$, which is the set of tempered distributions for which the integral norm
\begin{equation}
\label{eq:sobolev_p}
\| f\|_{W_s^p(\mathbb{R}^3)} =\|((1+|\,\cdot\,|^2)^\frac{s}{2} \widehat{f})\,\widecheck{ }\,\|_{L^p(\mathbb{R}^3)},
\end{equation}  
is finite. Here, the exponent $s$ controls the regularity of elements in $W_s^p(\mathbb{R}^3)$, and for $s>0$ multiplying the weight $(1+|\xi|^2)^\frac{s}{2}$ amounts to taking derivatives of the function $f$, while for $s<0$ it acts as smoothing on $f$. More precisely, the Sobolev space $W_s^p(\mathbb{R}^3)$ consists of all tempered distributions for which the $L^p$ norm is finite after convolution with the Bessel kernel $K_s$ defined by $\widehat{K_s}(\xi)=(1+|\xi|^2)^\frac{s}{2}$ \cite{Mattila2015}. Then, by the relation (\ref{eq:plancherel}), the norm of the Sobolev space for the case $p=2$, denoted by $H_{s}(\mathbb{R}^3)$, is equivalently defined as the integral on the frequency domain
\begin{equation}
\label{eq:sobolev_2}
\| f\|_{H_{s}(\mathbb{R}^3)} =\biggl( \int_{\mathbb{R}^3} (1+|\xi|^2)^s |\widehat{f}(\xi)|^2 d\xi\biggr)^\frac{1}{2},
\end{equation}
and it enables to compute the norm $\| f\|_{H_s(\mathbb{R}^3)}$ without the inversion formula (\ref{eq:ft_inverse}).
%
%
%
%
%
%
%
\subsection{Surface measures}
The present work is concerned with surface measures of two-dimensional (2D) surfaces in $\mathbb{R}^3$, which are also Borel and define integrals of measurable functions with respect to their area elements. For example, the surface measure of the $2$-sphere $\mathbb{S}^2 = \{x\in\mathbb{R}^3\,|\, |x|=1\}$, denoted by $\sigma_{\mathbb{S}^2}$, gives rise to integrals
\begin{equation}
\label{eq:def_surf_measure_sphere}
\int_{\mathbb{R}^3} f(x) d\sigma_{\mathbb{S}^2}(x) = \int_{\mathbb{S}^{2}}  f(x) dS,
\end{equation}
where $dS$ is the area element of $\mathbb{S}^2$, and its Fourier transform is known as
\begin{equation}
\label{eq:ft_nsphere}
\widehat{\sigma_{\mathbb{S}^2}}(\xi) = \int_{\mathbb{S}^{2}}  e^{-2\pi i \xi \cdot x} dS = \frac{2\pi}{|\xi|^{\frac{1}{2}}}J_{\frac{1}{2}}(2\pi|\xi|).
\end{equation}
Here, the function $J_\nu$ is the Bessel function of order $\nu$ whose series expansion is given in terms of the gamma function $\Gamma$:
\begin{equation}
\label{eq:Bessel_series}
J_\nu (t) = \sum^\infty_{j=0} \frac{(-1)^j}{j!} \frac{1}{\Gamma(j+\nu+1)}\biggl(\frac{t}{2} \biggr)^{2j+\nu}.
\end{equation}
Combining (\ref{eq:Bessel_series}) with $\Gamma(\frac{1}{2}+n) = \frac{(2n-1)!!}{2^n}\sqrt{\pi}$, the formula (\ref{eq:ft_nsphere}) is simplified \cite{Demeter2020} as
\begin{equation}
\label{eq:ft_nsphere_sinc}
\widehat{\sigma_{\mathbb{S}^2}}(\xi) =\frac{2\sin(2\pi|\xi|)}{|\xi|},
\end{equation} 
which claims that the Fourier transform $\widehat{\sigma_{\mathbb{S}^2}}$ is smooth on $\mathbb{R}^3$ but has the slow decay rate $\widehat{\sigma_{\mathbb{S}^2}}(\xi)=\mathcal{O} (|\xi|^{-1})\,(|\xi|\rightarrow \infty)$ as a singular measure.\par

More generally, the Paley-Wiener theorem \cite{Hormander2015} states that the Fourier transforms of finite Borel measures of compact support are smooth and extended to entire functions, which allows us to assume that compact surfaces in $\mathbb{R}^3$ have the associated surface measures whose Fourier transforms are smooth on the whole frequency domain. On the other hand, finiteness of the Sobolev norm of Borel measures in $\mathbb{R}^d$ is examined via the so-called $\alpha$-energy \cite{Mattila2015}: 
\begin{equation}
\label{eq:def_alpha_energy}
I_\alpha(\mu) = \iint \frac{1}{|x-y|^\alpha}d\mu(x)d\mu(y)= c(d,\alpha) \int_{\mathbb{R}^d}  |\xi|^{\alpha-d}|\widehat{\mu}(\xi)|^2d\xi,
\end{equation} 
where $0<\alpha<d$ and $c(d,\alpha)$ is a constant determined by $d$ and $\alpha$. Then, as a special class of the Hausdorff measures, the dimension of surface measures can be also understood as the supremum of the parameters $\alpha$ for which the integral (\ref{eq:def_alpha_energy}) is finite. Hence, it is clear that surface measures corresponding to 2D surfaces in $\mathbb{R}^3$ belong to the Sobolev space $H_s(\mathbb{R}^3)$ for $s<-\frac{1}{2}$. \par

Another relevant example is affine simplices, which are the main building blocks for modeling surfaces in applications. Let $Q_0^k$ be the standard $k$-simplex
\begin{equation}
\label{eq:simplex_st}
Q^k_0 = \{(x_1,\ldots,x_k)\in \mathbb{R}^k\,|\,x_1+\cdots +x_k \leq 1,\, x_i\geq 0 \,\,(i=1,\ldots,k)\},
\end{equation}
and an affine $k$-simplex in $\mathbb{R}^d\,(d\geq k)$ with vertices $\{p_0,p_1,\ldots, p_k\}$, denoted by $Q^k$, is the image of the standard simplex $Q_0^k$ by the affine map
\begin{equation}
\label{def_affine_simplex}
L(x) = p_0 + Ax,\quad x\in Q^k_0,
\end{equation}
where the $d\times k$ matrix $A$ satisfies $Ae_i = p_i - p_0\,(i=1,\ldots,k)$. Then, the Fourier transforms of the surface measures of affine simplices also have closed forms for general $d$ and $k$ \cite{JiWaHuMaNi2019}, and the present work is in particular related to the case for $d=3$ and $k=2$. By an abuse of notation, the surface measure of a $k$-simplex $Q^k$ is denoted by the same symbol, and its Fourier transform is written as
\begin{equation}
\label{eq:simplex32_ft}
\widehat{Q^k}(\xi) = -\frac{\gamma}{4\pi^2} \sum^2_{i=0}\frac{e^{-2\pi i\xi\cdot p_i}}{\prod_{j\neq i}\{\xi\cdot (p_i-p_j)\} },\quad \xi \in \mathbb{R}^3,
\end{equation}
where $\gamma$ is the ratio of the surface area of $Q^k$ to that of $Q^k_0$. Again, the function (\ref{eq:simplex32_ft}) is smooth due to the compactness of $Q^k$, and the corresponding measure belongs to the space $H_s(\mathbb{R}^3)$ for $s<-\frac{1}{2}$.
%
%
%
%
%
%
%
\subsection{Objective function and derivatives}
Now suppose that we are given two surface measures $\sigma_1$ and $\sigma_2$ corresponding to compact surfaces $S_1$ and  $S_2$ in $\mathbb{R}^3$  and consider local minimization of the distance
\begin{align}
\label{eq:def_distance_Sobolev}
d_{H_s}(\sigma_1, \sigma_{2,\varphi(X)}) =\|\sigma_1 - \sigma_{2,\varphi(X)} \|_{H_{s}(\mathbb{R}^3)},
\end{align}
for appropriate $s<-\frac{1}{2}$ with respect to the vector $X$ via the objective function
\begin{align}
\label{eq:def_objective}
f(X) =  \int_{\mathbb{R}^3} (1+|\xi|^2)^s |\widehat{\sigma_1}(\xi)-\widehat{\,\sigma_{2,\varphi(X)\,}}(\xi)|^2 d\xi.
\end{align}
Here, the measure $\sigma_{2,\varphi(X)}$ is the pushforward of $\sigma_2$ by the isometry $\varphi(X)$ that lies in the special Euclidean group $SE(3)$. In particular,  by setting $X = (b,Y)$ where 
\begin{equation}
\label{eq:def_skew}
b=
\begin{pmatrix}
b_1 \\
b_2 \\
b_3\\
\end{pmatrix},\quad
Y=
\begin{pmatrix}
0 & y_1 & y_2 \\
-y_1 & 0 & y_3 \\
-y_2 & -y_3 &0\\
\end{pmatrix},
\end{equation}
elements of $SE(3)$ as a Lie group can be uniquely represented as
\begin{equation}
\label{eq:def_rigid_close1}
[\varphi(X)](x) = b+e^Yx,\quad x \in \mathbb{R}^3,
\end{equation}
in small neighborhoods of the identity element. Here, the exponential map $e^Y$ serves as a local diffeomorphism from the space of skew-symmetric matrices to the special orthogonal group $SO(3)$ \cite{Stillwell2008}, and allows to translate the problem on the Lie group into that on its Lie algebra without causing large distortions around the identity.\par

Next, we derive a formula for the gradient $\nabla f$ with respect to the variables $X=(b,Y)$. For each $x_i$ in $X$, the partial derivative $\partial f / \partial x_i$ is generally written as 
\begin{align}
\label{eq:deriv_obj_1st}
\nonumber   \frac{\partial f}{\partial  x_i} &=\frac{\partial}{\partial x_i}\int_{\mathbb{R}^3} \widehat{K_s}^2(\xi) |\widehat{\sigma_1}(\xi)-\widehat{\sigma_{2,\varphi(X)}}(\xi)|^2 d\xi\\
\nonumber                            &=\frac{\partial}{\partial x_i }\int_{\mathbb{R}^3} \widehat{K_s}^2(\xi)  (\widehat{\sigma_1}(\xi)-\widehat{\sigma_{2,\varphi(X)}}(\xi))(\overline{\widehat{\sigma_1}(\xi)-\widehat{\sigma_{2,\varphi(X)}}(\xi)}) d\xi\\
\nonumber                            &=\int_{\mathbb{R}^3} \widehat{K_s}^2(\xi)  \biggl({\widehat{\sigma_1}(\xi)
                                                  -\frac{\partial \widehat{\sigma_{2,\varphi(X)}}}{\partial x_i}(\xi)}\biggr)(\overline{\widehat{\sigma_1}(\xi)-\widehat{\sigma_{2,\varphi(X)}}(\xi)}) d\xi\\
\nonumber                            &+\int_{\mathbb{R}^3} \widehat{K_s}^2(\xi) (\widehat{\sigma_1}(\xi)-\widehat{\sigma_{2,\varphi(X)}}(\xi))\biggl(\overline{\widehat{\sigma_1}(\xi)
                                                  -\frac{\partial \widehat{\sigma_{2,\varphi(X)}}}{\partial x_i}(\xi)}\biggr)d\xi\\
                                             &=2\int_{\mathbb{R}^3} \widehat{K_s}^2(\xi) \cdot \mathrm{Re}\biggl\{(\widehat{\sigma_1}(\xi)-\widehat{\sigma_{2,\varphi(X)}}(\xi))\biggl(\overline{\widehat{\sigma_1}(\xi)
                                                  -\frac{\partial \widehat{\sigma_{2,\varphi(X)}}}{\partial x_i}(\xi)}\biggr)\biggr\} d\xi,
\end{align}
where $\widehat{K_s}(\xi) = ( 1+ |\xi|^2)^\frac{s}{2}$. Hence, in order for the formula (\ref{eq:deriv_obj_1st}) to make sense, it is also necessary to compute $x_i$-derivatives of the function 
\begin{align}
\label{eq:fourier_under_phi}
\widehat{\,\sigma_{2,\varphi(X)}\,}(\xi) = \int_{\mathbb{R}^3}  e^{-2\pi i \xi \cdot [\varphi(X)](x)}d\sigma_2(x)= e^{-2\pi i \xi \cdot b}\int_{\mathbb{R}^3}  e^{-2\pi i \xi \cdot (e^Yx)} d\sigma_2(x).
\end{align}
For the translation $b$, one can easily obtain
\begin{align}
\label{eq:deriv_fourier_under_phi_bi}
\frac{\partial \widehat{\sigma_{2,\varphi(X)}}}{\partial b_i}(\xi)  = -2\pi i \xi_i e^{-2\pi i \xi \cdot b}\int_{\mathbb{R}^3}  e^{-2\pi i \xi \cdot (e^Y x)} d\sigma_2(x),
\end{align}
while differentiating it with respect to $y_i$ leads to
\begin{align}
\label{eq:deriv_fourier_under_phi_yi}
\nonumber \frac{\partial \widehat{\sigma_{2,\varphi(X)}}}{\partial y_i}(\xi) &= e^{C\xi \cdot b}\int_{\mathbb{R}^3}\biggl\{C\xi \cdot \biggl(\frac{\partial e^Y}{\partial y_i}x\biggr)\biggr\}  e^{C\xi \cdot (e^Yx)} d\sigma_2(x)\\
\nonumber                                                                                    &= Ce^{ C\xi \cdot b}\int_{\mathbb{R}^3}\biggl[ \xi \cdot \biggl\{\biggl(\frac{\partial Y}{\partial y_i}e^Y\biggr)x\biggr\}\biggr]  e^{C\xi \cdot (e^Yx)} d\sigma_2(x)\\
\nonumber                                                                                   &= Ce^{ C\xi \cdot b}\int_{\mathbb{R}^3}\biggl[ \biggl\{\biggl(\frac{\partial Y}{\partial y_i}e^Y\biggr)^{\mathsf{T}} \xi\biggr\} \cdot x\biggr]  e^{C\xi \cdot (e^Yx)} d\sigma_2(x)\\
		                                                                                    &= C e^{C\xi \cdot b}\sum^3_{k=1}F_{i,k}(Y, \xi) \int_{\mathbb{R}^3}x_k e^{C \xi \cdot (e^Yx)} d\sigma_2(x),
\end{align}
where $C =-2\pi i$ and the factors $F_{i,k}$ are independent of the variable $x$. Thus, to evaluate the gradient $\nabla f$ for $Y$, it suffices to have additional three transforms
\begin{align}
\label{eq:def_fourier_weight}
 \int_{\mathbb{R}^3}x_k e^{-2\pi i \xi \cdot (e^Yx)} d\sigma_2(x),\quad \xi \in \mathbb{R}^3\quad (k=1,2,3)
\end{align}
which are essentially the Fourier transforms with the monomial weights $x_k$. Likewise, the second derivatives can also be found in terms of the Fourier-type transforms
\begin{align}
\label{eq:fourier_yiyj}
\int_{\mathbb{R}^3}x_k x_l e^{C \xi \cdot (e^Yx)} d\sigma_2(x),\quad (1\leq k \leq l \leq 3) 
\end{align}
although we avoid directly evaluating the Hessian $\nabla^2 f$ and resort to employing an optimization algorithm that works with the gradient information only.


\section{Numerical methods}\label{method}
This section is devoted to describing numerical methods for computing the objective function (\ref{eq:def_objective}) with triangulated surfaces and locally optimizing it on the group $SE(3)$ via the representation (\ref{eq:def_rigid_close1}). The former consists of numerical quadratures on the spatial and frequency domains, which are bridged by a fast summation algorithm called NUFFT, while the latter is based on the trust-region method with the symmetric rank-1 (SR1) Hessian updates.
%
%
%
%
%
%
%
\subsection{Numerical quadratures on simplex}
As shown in Section \ref{formulate}, the most essential step in local optimization of the objective function (\ref{eq:def_objective}) is to approximate the Fourier transforms $\widehat{\sigma_1}$ and $\widehat{\,\sigma_{2,\varphi(X)\,}}$, respectively, and the three Fourier-type transforms (\ref{eq:def_fourier_weight}) in the derivatives (\ref{eq:deriv_fourier_under_phi_yi}). Due to the linearity, this task for triangulated surfaces is reduced to computing the Fourier transform of each simplex and summing up their contributions on the frequency domain. While the exact formula (\ref{eq:simplex32_ft}) and its derivatives with respect to vertex positions \cite{IcSiDeTu2025} provide the most direct way for the purpose, such closed forms may suffer from large cancellation errors and limit accuracy in practical computations \cite{JiWaHuMaNi2019, IcSiDeTu2025}. \par

Following Koga \cite{Koga2024}, the present work chooses discretization of affine simplices by fast-converging numerical quadratures on the standard simplex (\ref{eq:simplex_st}). Namely, with the affine map $L$ given by (\ref{def_affine_simplex}), we perform the change of variables
\begin{equation}
\label{eq:change_variable}
\int_{Q^k} f(x) dS = \int_{Q_0^k} f(L(u)) \gamma_A du,
\end{equation}
where the constant $\gamma_A$ is the Jacobian determined by the matrix $A$, and a numerical quadrature $\{(u_i,q_i)\}$ on $Q^k_0$ is applied to the right-hand side of (\ref{eq:change_variable}): 
\begin{equation}
\label{eq:change_variable_approx}
 \int_{Q_0^k} f(L(u)) \gamma_A du \approx \sum_{i}  f(L(u_i)) q_i \cdot \gamma_A.
\end{equation}
To achieve nearly optimal approximations, our implementation employs the class of symmetric rules by Xiao and Gimbutas \cite{XiGi2010} on the regular triangle with the vertices $(-1,-1/\sqrt{3})$, $(0,2/\sqrt{3})$, and $(1,-1/\sqrt{3})$ and adapts those as $\{(u_i,q_i)\}$ on the standard domain $Q^k_0$, while the iterative use of the Gauss-Legendre quadratures (i.e., tensor product rules) also works with suboptimal efficiency.\par

Here, we remark that the nodes $\{L(u_i)\}$ and weights $\{q_i \cdot \gamma_A\}$ in (\ref{eq:change_variable_approx}) are independent of the integrand $f$. In fact, this procedure is equivalent to converting the numerical quadrature $\{(u_i,q_i)\}$ on $Q^k_0$ to $\{(L(u_i),q_i \cdot \gamma_A)\}$ on $Q^k$, and we discretize the whole surface measure of a triangulated surface in terms of the Dirac measures on $\mathbb{R}^3$ with nonuniform coefficients. Although such a weighted point cloud is no longer geometric, we see later that this viewpoint plays a crucial role in efficiently summing up their contributions on the frequency domain.
%
%
%
%
%
%
%
\subsection{Numerical quadratures on $\mathbb{R}^3$}
Next, we select a numerical quadrature on the whole $\mathbb{R}^3$ that approximates the objective function (\ref{eq:def_objective}) and its derivatives (\ref{eq:deriv_obj_1st}). Again, following \cite{Koga2024}, let $D=[-\xi_\text{max},\xi_\text{max}]^3$ and consider the lattice
\begin{equation}
\label{eq:def_trap_grid}
\xi_{i,j,k} = (ih,jh,kh),\quad h=\frac{2\xi_\text{max}}{N}, \quad i,j,k=-\frac{N}{2},\ldots, \frac{N}{2}.
\end{equation}
Then, the 3D trapezoidal rule is applied to functions on $D$ as
\begin{equation}
\label{eq:trap_rule_R3}
\int_D f(\xi) d\xi\approx  \sum^{\frac{N}{2}}_{i=-\frac{N}{2}}\sum^{\frac{N}{2}}_{j=-\frac{N}{2}}\sum^{\frac{N}{2}}_{k=-\frac{N}{2}} f(\xi_{i,j,k}) w_{i,j,k} h^3,
\end{equation}
where $w_{i,j,k}$ is 1 for points in the interior of the cube $D$, $1/2$ for the faces, $1/4$ for the edges, and $1/8$ for the corners.\par

To understand errors from the formula (\ref{eq:trap_rule_R3}), we refer to the exactness of the trapezoidal rules in relation to the Fourier transforms \cite{TrWe2014}. That is, for well-behaved integrands $f$ on $\mathbb{R}^3$, let $I_h(f)$ be the infinite trapezoidal sum of width $h$
\begin{equation}
\label{eq:def_trap_sum}
I_h(f)  = h^3\sum_{k\in \mathbb{Z}^3} {f} (hk),
\end{equation}
and its discretization errors can be written in terms of the Fourier transform $\widehat{f}$:
\begin{equation}
\label{eq:error_trap}
I_h(f) - I(f) =  \sum_{k\in \mathbb{Z}^3\backslash \{0\}} \widehat{f} (k/h),
\end{equation}
where $I(f)$ is the integral of $f$ and is equal to the value $\widehat{f}(0)$. This is an immediate consequence of the Poisson summation formula \cite{Grafakos2014-1, Henrici1977}, and it states that the sum (\ref{eq:def_trap_sum}) serves as an exact numerical quadrature for functions whose Fourier transforms are supported on compact sets within $(-1/h, 1/h)^3$. In the present work, when switching the roles of $x$ and $\xi$, surface measures are convoluted with the Bessel kernel and the resulting functions on the spatial domain are no longer compactly supported. Nevertheless, by properly choosing $h$ and $\xi_\text{max}$ that both depend on the regularity $s$, the finite sum (\ref{eq:trap_rule_R3}) yields good approximations of the integrals (\ref{eq:def_objective}) and (\ref{eq:deriv_obj_1st}). 
%
%
%
%
%
%
%
\subsection{Nonuniform fast Fourier transforms}
To obtain data on the lattice (\ref{eq:def_trap_grid}) from a weighted point cloud $\{(x_l, c_l) \}^M_{l=1}$, we evaluate finite sums in the form
\begin{alignat}{2}
\label{eq:def_sum_se3}
\sum_{l=1}^{M}c_le^{-2\pi i\xi \cdot (b+e^Yx_l)},
\end{alignat}
where the matrix exponential $e^Y$ is computed via the scaling and squaring method \cite{Higham2009}. Most naively, the direct summation of contributions from $M$ source points has $\mathcal{O}(M\cdot N^3)$ complexity, and the fast Fourier transform (FFT) is not straightforwardly applicable to nonuniform data in the spatial domain. Instead, an efficient algorithm called the Type-1 nonuniform fast Fourier transform (Type-1 NUFFT) is capable of approximating the following sum to a prescribed relative accuracy $\epsilon_\text{rel}$: 
\begin{alignat}{2}
\label{eq:def_nufft_type1}\widehat{f}_\eta=\sum_{l=1}^{M}f_le^{-i\eta \cdot x_l}, \quad & \eta=(\eta_1, \eta_2, \eta_3)\in \mathbb{Z}^3 \cap [-\eta_\text{max}, \eta_\text{max}]^3,
\end{alignat}
where the coefficients $\{f_l\}$ are any complex numbers and the coordinates $\{x_l\}$ lie in $[-\pi,\pi)^3$. To utilize FFT, the Type-1 NUFFT regards the sum (\ref{eq:def_nufft_type1}) as the exact Fourier transform of atomic measures $ \mu = \sum_{l=1}^{M}f_l\delta_{x_l}$ and replaces it with
\begin{equation}
\label{def_nufft_conv}
 \mu * \psi = \sum_{l=1}^{M}f_l \psi(\cdot-x_l),
\end{equation}
where $*$ is the convolution on the torus $\mathbb{T}^3=[-\pi,\pi)^3$ and $\psi$ is a smooth periodic function with a small numerical support. Then, the Fourier transform of the regularized measure is computed by the 3D trapezoidal rule with the FFT acceleration, and dividing the results by the symbols $\widehat{\psi}(\eta)$ yields approximations of the desired sums. Thus, as opposed to the direct summation, the complexity of this procedure is $\mathcal{O}(N^3\log N^3 + M \log^3 (1/\epsilon_\text{rel}))$. In the present work, we employ the CPU version of FINUFFT \cite{BaMaKl2019}, which chooses $\psi$ to be the ``exponential of semicircle" kernel with nearly optimal aliasing errors \cite{Barnett2021}. Although the trapezoidal rule on $\mathbb{R}^3$ requires integrand values on $h\mathbb{Z}^3$ rather than $\mathbb{Z}^3$, by rewriting the sum of interest as
\begin{alignat}{2}
\label{eq:source_scale}
\sum_{l=1}^{M}f_le^{-ih\eta \cdot x_l}= \sum_{l=1}^{M}f_le^{-i\eta \cdot  hx_l},
\end{alignat}
it turns out that the evaluations on $h\mathbb{Z}^3$ are reduced to NUFFT with the coordinates $\{hx_l\}$ \cite{Koga2024}. In other words, the Fourier transform on the lattice (\ref{eq:def_trap_grid}) for $0<h<1$ is obtained from a more clustered point cloud in the vicinity of the origin.
%
%
%
%
%
%
%

\subsection{SR1 trust-region method}
We locally minimize the objective function (\ref{eq:def_objective}) using the trust-region method with the symmetric rank-1 (SR1) Hessian updates \cite{NoWr2006}. Namely, the solution $x_k$ at the step $k$ is updated to $x_{k+1} = x_k+p$ by finding the minimizer $p$ of a quadratic model around $x_k$ in the form
\begin{alignat}{2}
\label{eq:trust_subproblem}
m_k(p) = f(x_k)+ g_k^{\mathsf{T}} p+p^{\mathsf{T}}B_{k} p,\quad \|p\|_2 \leq \Delta_k.
\end{alignat}
Here, $\|\cdot\|_2$ is the Euclidean norm on $\mathbb{R}^n$, $\Delta_k$ is the trust-region radius, $g_k= \nabla f(x_k)$, and  $B_k$ is a symmetric matrix as an approximation of the true Hessian $\nabla^2 f(x_k)$. More specifically, we begin with $B_0 = I_n$ and improve $B_k$ by the SR1 update formula 
\begin{alignat}{2}
\label{eq:sr1_update}
B_{k+1} = B_k + \frac{(y_k - B_k s_k)(y_k - B_k s_k)^{\mathsf{T}}} { s_k \cdot (y_k - B_k s_k)},
\end{alignat}
where $s_k = x_{k+1}-x_k$ and $y_k = \nabla f(x_{k+1})-\nabla f(x_{k})$, if the condition
\begin{alignat}{2}
\label{eq:sr1_condition}
\frac{s_k \cdot (y_k - B_k s_k)}{\|s_k\|_2 \|y_k - B_ks_k\|_2} > 10^{-8},
\end{alignat}
holds, and we set $B_{k+1}=B_k$ otherwise. In general, the SR1 update does not guarantees that $B_{k+1}$ is positive definite even if $B_k$ is so, and hence it naturally leads us to the trust-region framework that involves a few matrix decompositions. Nevertheless, the dimension of the optimization problem is rather small and solving the subproblem (\ref{eq:trust_subproblem}) is inexpensive in comparison to evaluating the Fourier transform $\widehat{\,\sigma_{2,\varphi(X)\,}}$ and its derivatives. Besides, the formula (\ref{eq:sr1_update}) does not assume the curvature condition that requires additional computations of the gradient $\nabla f$, and only values $f(x_k+p)$ are needed in quality estimations of the models $m_k$. To update the radius $\Delta_k$, we follow the SR1 trust-region algorithm described in \cite{NoWr2006}, where the priority is given to improvements on $B_k$ rather than substantial decreases in the function $f$.  \par 

At each step, we first perform the eigendecomposition of the matrix $B_k$ and compute the solution $p = -B^{-1}_k g_k$ if its eigenvalues are all positive. Then, if $B_k$ is not positive definite or the constraint $\|p\|_2\leq\Delta_k$ is violated, the minimizer $p$ on the boundary $\|p\|_2=\Delta_k$ is sought using the stable algorithm suggested by Adachi, Iwata, Nakatsukasa, and Takeda \cite{AdIwNaTa2017}, which replaces the classical iterative processes \cite{NoWr2006} with solving just one generalized eigenvalue problem in the form
\begin{alignat}{2}
\label{eq:def_general_eigen}
\begin{pmatrix}
-I_n & B_k  \\
B_k & -g g^{\mathsf{T}}/ \Delta^2_k
\end{pmatrix}
u=-\lambda
\begin{pmatrix}
0 & I_n  \\
I_n & 0
\end{pmatrix}
u.
\end{alignat}
Their analysis shows that the right-most eigenvalue $\lambda^*$ is in fact the Lagrange multiplier for the desired solution and, if the subproblem is not ``hard", the vector $p$ can be constructed from the eigenvector corresponding to $\lambda^*$. On the other hand, in the hard case, we exploit data from the eigendecompostion of $B_k$ and obtain $p$ in the standard way \cite{NoWr2006}. In our codes, these operations are implemented with \texttt{Eigen} \cite{eigenweb}. \par

The strategy above is somewhat different from one in \cite{AdIwNaTa2017}, where the equations $B_kp = - g_k$ are solved without checking the eigenvalues of $B_k$ and the hard case is treated with data from the generalized eigenvalue problem (\ref{eq:def_general_eigen}). However, performing such decompositions for at most $n=12$ has little impact on overall computational costs, while it prevents us from applying the conjugate gradient method to general symmetric matrices and also gives a simpler implementation for the hard case. 

\section{Numerical results}\label{results}
Now we show numerical results to validate the algorithm described in Section \ref{method}. Throughout this section, one triangulated surface is placed in $\mathbb{R}^3$ while its copy is rigidly transformed by an isometry $\varphi(X)$ in the form (\ref{eq:def_rigid_close1}), where the former is denoted by $S_1$ and the latter is $S_2$, and we ask whether one can find the inverse $\varphi(X)^{-1}$ as a trivial zero of the objective function (\ref{eq:def_objective}) via numerical optimization. Without loss of generality, we can set the barycenter of the fixed surface $S_1$ to the origin before duplicating it as $S_2$. It is not only because relative positioning between the given objects is essential, but also because accuracy of the trapezoidal rule on the frequency domain is dictated by how far the corresponding functions on the spatial domain spread from the origin. In the following, the implemented codes in C++ with OpenMP are run on MacBook Pro 2024 (Apple M4 Max and 36 GB memory). As parameters, we fix $N=64$ and $\xi_\text{max} = 5$ for (\ref{eq:trap_rule_R3}), the relative tolerance for the Type-1 NUFFT algorithm is set to $\epsilon_\text{rel}=10^{-15}$, and the optimization algorithm terminates if the condition $\|\nabla f\|_2 < 10^{-7}$ holds, unless otherwise specified.
%
%
%
%
%
%
%
\subsection{Symmetric case}
\begin{figure}[t]
\begin{center}
 \includegraphics[width=0.8\linewidth]{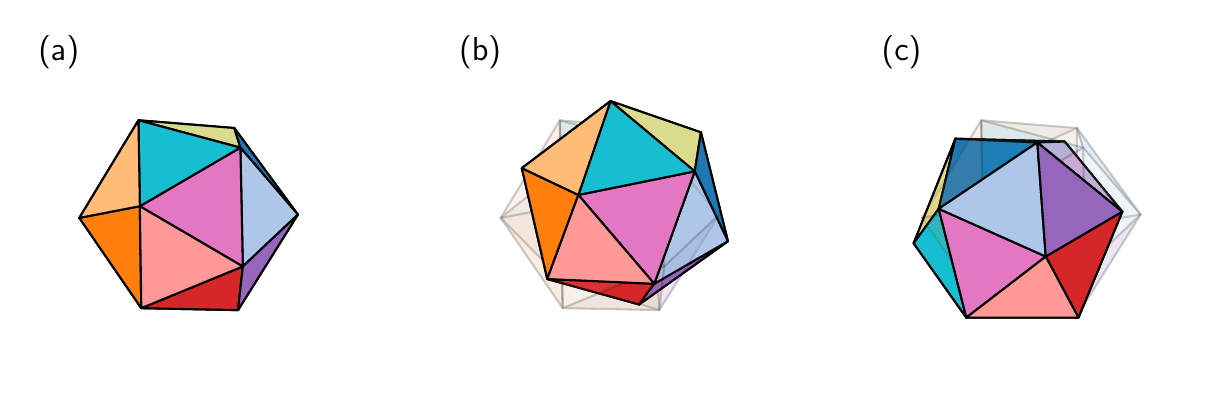}
\end{center}
\caption{\label{fig:plot_icosad_three} Illustrations of triangulated surfaces in case of two icosahedra: (a) fixed surface $S_1$,  (b) initial position of $S_2$ from (\ref{eq:def_bY_icosad1}), and (c) that from (\ref{eq:def_bY_icosad2}).} 
\end{figure}

We start numerical experiments with the case of two regular icosahedra for $s=-10$. First, we pick the parameters $X=(b,Y)$ to be
\begin{equation}
\label{eq:def_bY_icosad1}
b_1=b_2=b_3 = 0.1,\quad y_1 = \frac{\pi}{16},\,\,y_2 = \frac{\pi}{9}, \,\,y_3 =0.
\end{equation}
Here, Figure \ref{fig:plot_icosad_three}(a) shows as $S_1$ a regular icosahedron inscribed to the unit sphere $\mathbb{S}^2$. On the other hand, we apply the isometry $\varphi(X)$ from (\ref{eq:def_bY_icosad1}) to a copy of $S_1$ and treat it as $S_2$, which is illustrated in Fig.\ref{fig:plot_icosad_three}(b) with the transparent $S_1$. Moreover, we put unique colors on the faces of both surfaces so that the correspondence between the vertices are clear. As one can see, the two objects are visually distinct as sets in $\mathbb{R}^3$, whereas they are kept sufficiently close so that the suggested algorithm can capture the inverse $\varphi(X)^{-1}$. For these rough surfaces, each triangle of $S_1$ is discretized using the 171-point rule on (\ref{eq:simplex_st}), and the 79-point rule is applied to $S_2$. \par

\begin{figure}[t]
\begin{center}
 \includegraphics[width=\linewidth]{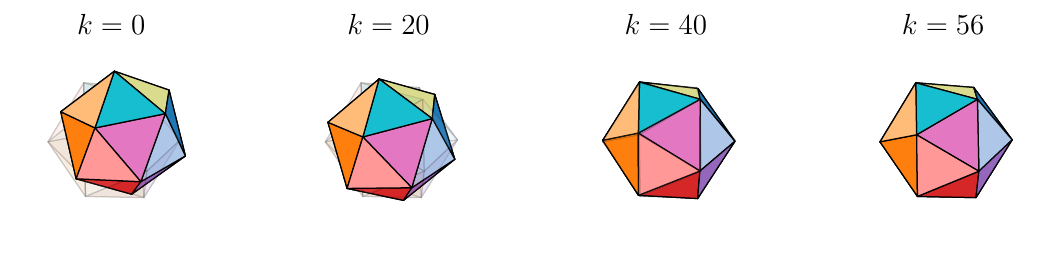}
\end{center}
\caption{\label{fig:snaps_icosad_inverse} Snapshots of surface $S_2$ in case of two icosahedra at indicated steps for example (\ref{eq:def_bY_icosad1}).} 
\end{figure}

To see how the surface $S_2$ behaves during optimization, Figure \ref{fig:snaps_icosad_inverse} shows snapshots of $S_2$ at a few steps $k$. In this case, it seems that the translation part $b$ is optimized first and the inverse $\varphi(X)^{-1}$ is successfully obtained by subsequent minimization on the rotation part $Y$. We speculate that it is because the translation plays a major role if the barycenters of the two surfaces are far enough while the rotation is more significant once the objects become sufficiently close. Although the present work prioritizes generality of the problem, it is in fact an option to omit the variable $b$ from the formulation in global alignment \cite{SiYa2024}, where it often suffices to minimize an objective function on the group $SO(3)$ after setting the barycenters to the origin. However, the situation is more complicated in partial alignment \cite{RiZhCoChDu2025} and generally the translation part cannot be dropped there.\newpage

\begin{figure}[t]
\begin{center}
 \includegraphics[width=0.8\linewidth]{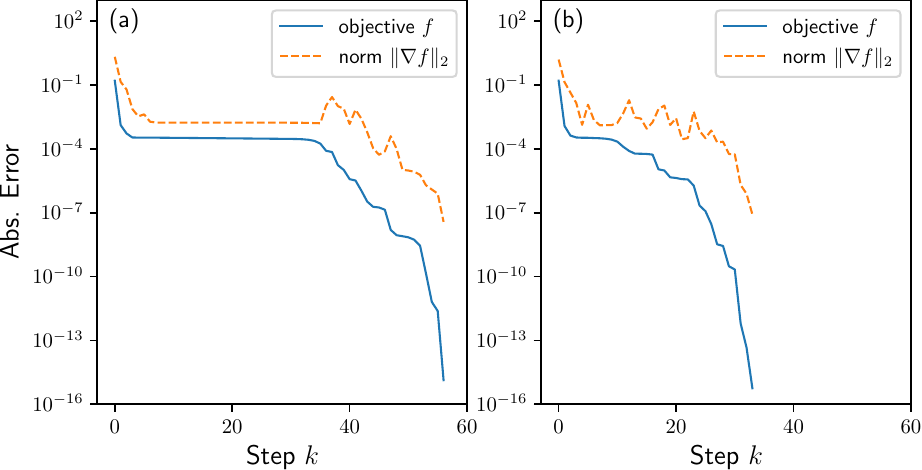}
\end{center}
\caption{\label{fig:conv_icosad} Convergence of SR1 trust-region method in case of two icosahedra for $s=-10$: (a) objective function and gradient norm for initial position from (\ref{eq:def_bY_icosad1}) and (b) those from (\ref{eq:def_bY_icosad2}).} 
\end{figure}

Next, we take a close look at behavior of the objective function and its gradient norm. For this purpose, Figure \ref{fig:conv_icosad}(a) shows absolute values of $f$ (solid line) and $\|\nabla f\|_2$ (dashed line) versus the step $k$. As easily seen, the two curves have clear negative slopes at the first few steps and soon become almost flat until $k\approx 35$, while the monotonicity of $f$ is guaranteed as a property of the trust-region method. Then, once the solution $(b_k, Y_k)$ approaches a root, substantial decreases in both quantities are observed again and continue until the norm $\|\nabla f\|_2$ reaches the termination condition at a seemingly superlinear rate.

\begin{figure}[t]
\begin{center}
 \includegraphics[width=\linewidth]{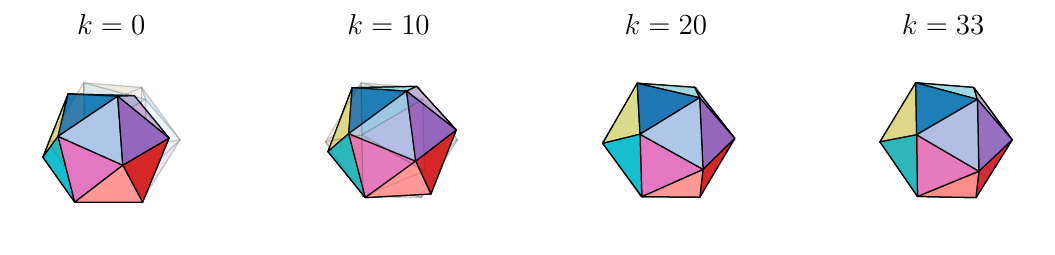}
\end{center}
\caption{\label{fig:snaps_icosad_nonconvex} Snapshots of surface $S_2$ in case of two icosahedra at indicated steps for example (\ref{eq:def_bY_icosad2}).} 
\end{figure}

To understand the objective function (\ref{eq:def_objective}) in more detail, we also consider
\begin{equation}
\label{eq:def_bY_icosad2}
b_1=b_2=b_3 = -0.1,\quad y_1 = 0,\,\,y_2 = \pi+\frac{\pi}{8}, \,\,y_3 =\pi+\frac{\pi}{16},
\end{equation}
and the corresponding $S_2$ is shown in Fig.\ref{fig:plot_icosad_three}(c). In this case, the rotated surface is so far from its original position that derivative-based optimization should fail to seek the inverse $\varphi(X)^{-1}$ due to the icosahedral symmetry. Again, Figure \ref{fig:snaps_icosad_nonconvex} shows a few snapshots of the surface $S_2$, and optimization in the translational direction is performed first. As opposed to the case (\ref{eq:def_bY_icosad1}), however, 
 the suggested algorithm returns the identical position as a set in $\mathbb{R}^3$ but the colors of the faces are permuted. In fact, Figure \ref{fig:conv_icosad}(b) shows behavior of the function $f$ and its gradient norm for the case (\ref{eq:def_bY_icosad2}), which is similar to that in Fig.\ref{fig:conv_icosad}(a) except for its much shorter flat region, and implies that the resulting isometry is another zero of the objective function but not the inverse of $\varphi(X)$ generated from (\ref{eq:def_bY_icosad2}). In other words, we can claim that the algorithm succeeds in finding the inverse $\varphi(X)^{-1}$ up to symmetry, and there are $n$ solutions in this example where $n$ is the order of the icosahedral group. Thus, it is concluded that the optimization problem with the objective function (\ref{eq:def_objective}) may allow multiple global solutions and hence is nonconvex in general, and the property is easily confirmed using surfaces with nontrivial symmetry.\par

\begin{figure}[t]
\begin{center}
 \includegraphics[width=0.8\linewidth]{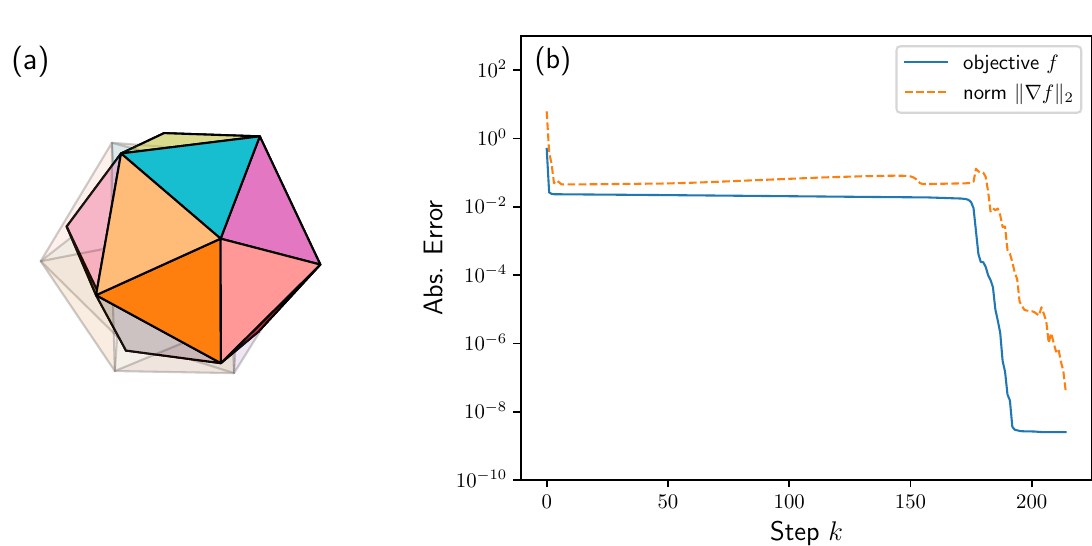}
\end{center}
\caption{\label{fig:conv_icosad_farthest} Hard problem with case of two icosahedra: (a) initial position of surface $S_2$ from (\ref{eq:def_bY_embryo}) and (b) objective function and gradient norm for $s=-5$.} 
\end{figure}

\begin{figure}[t]
\begin{center}
 \includegraphics[width=0.8\linewidth]{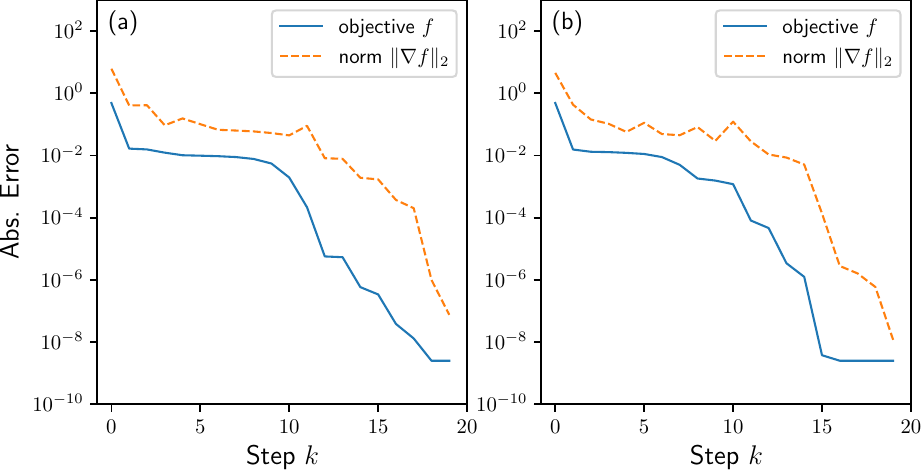}
\end{center}
\caption{\label{fig:conv_icosad_sharp} Convergence of SR1 trust-region method in case of two icosahedra for $s=-5$: (a) objective function and gradient norm for initial position from (\ref{eq:def_bY_icosad1}) and (b) those from (\ref{eq:def_bY_icosad2}).} 
\end{figure}

As shown above, the suggested algorithm may get stuck at points away from desired solutions and such phenomena can be considered the main cause of inefficiency. In our experience, the number of steps that is required for convergence is largely affected by how far the initial position is from stationary points, as in many other optimization problems. For instance, when the isometry $\varphi(X)$ from 
\begin{equation}
\label{eq:def_bY_embryo}
b_1=b_2=b_3 = 0.1,\quad y_1 = -\frac{\pi}{4},\,\,y_2 = 0, \,\,y_3 =0,
\end{equation}
is applied to the case of icosahedra, it yields the surface $S_2$ that is almost the farthest from $S_1$ in the rotational direction, as shown in Fig.\ref{fig:conv_icosad_farthest}(a), and the trust-region method fails to converge for $s=-10$ even after 1,000 steps. To improve the convergence property in practice, a promising option is to increase the regularity $s$, which enables the Sobolev norm to compare objects more sharply. In fact, Figure \ref{fig:conv_icosad_farthest}(b) shows the objective function $f$ and its gradient norm $\|\nabla f\|_2$ for the case of icosahedra with (\ref{eq:def_bY_embryo}) and $s=-5$ and indicates that the algorithm barely converges after getting stuck for many $k$. Moreover, the same pair of functions from the previous examples (\ref{eq:def_bY_icosad1}) and (\ref{eq:def_bY_icosad2}) is plotted in Fig.\ref{fig:conv_icosad_sharp} for $s=-5$, where one can easily see that convergence in both cases is significantly improved. Thus, efficiency of the suggested algorithm largely depends on the choice of the regularity $s$, and poor initial guesses may require large $s$ for convergence. However, we also remark here that implementing the algorithm becomes highly delicate as the parameter $s$ approaches its critical value that justifies the second-order optimization framework.

%
%
%
%
%
%
%
\subsection{Nonsymmetric case}
\begin{figure}[t]
\begin{center}
 \includegraphics[width=0.7\linewidth]{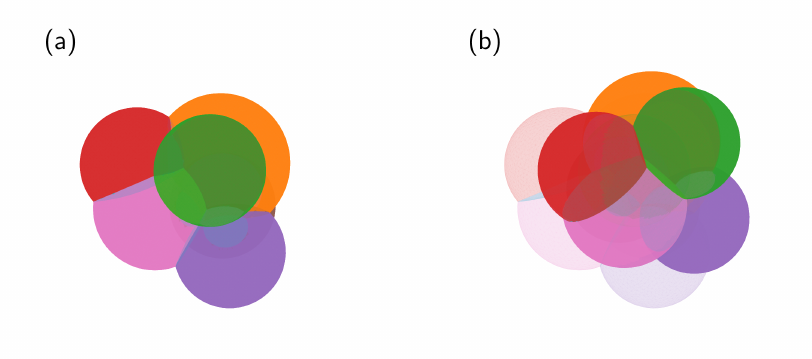}
\end{center}
\caption{\label{fig:plot_foam_two} Illustrations of triangulated surfaces in case of embryo model: (a) fixed surface $S_1$ and (b) initial position of $S_2$ from (\ref{eq:def_bY_embryo}).} 
\end{figure}

We proceed to a more practical example that is nonsymmetric and therefore gives a unique zero of the objective function. Figure \ref{fig:plot_foam_two}(a) shows as $S_1$ a foam-like surface with multiple junctions, and we put different colors on smooth components for clear illustration. This model consists of 23,055 triangles with 11,357 nodes and is provided as part of the \texttt{foambryo} package \cite{IcDeMcDuTu2023}, an end-to-end computational method for performing force inference on embryonic tissues. For this example, we select the parameters (\ref{eq:def_bY_embryo}) again and the resulting surface $S_2$ is shown in Fig.\ref{fig:plot_foam_two}(b). Unlike the case of the regular icosahedra, this embryo model has no symmetry and the inverse $\varphi(X)^{-1}$ is expected to be the unique minimizer. Moreover, the maximum diameter of triangles that constitute the surface is far smaller than the previous case. Utilizing this feature, we choose the 6-point rule for discretizing the moving surface $S_2$, while the 55-point rule is applied once to the fixed $S_1$. Such imbalanced sampling plays a crucial role in speeding up the overall computations and also emphasizes that comparisons are done as surfaces rather than as point clouds. \par

\begin{figure}[t]
\begin{center}
 \includegraphics[width=\linewidth]{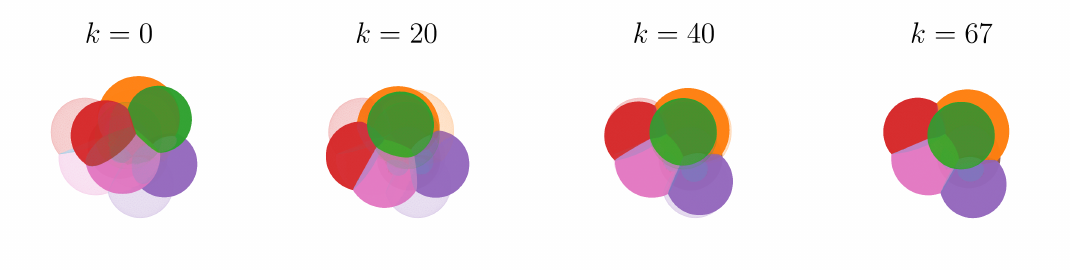}
\end{center}
\caption{\label{fig:snaps_foam} Snapshots of surface $S_2$ in case of embryo model at indicated steps for example (\ref{eq:def_bY_embryo}).} 
\end{figure}

\begin{figure}[t]
\begin{center}
 \includegraphics[width=0.8\linewidth]{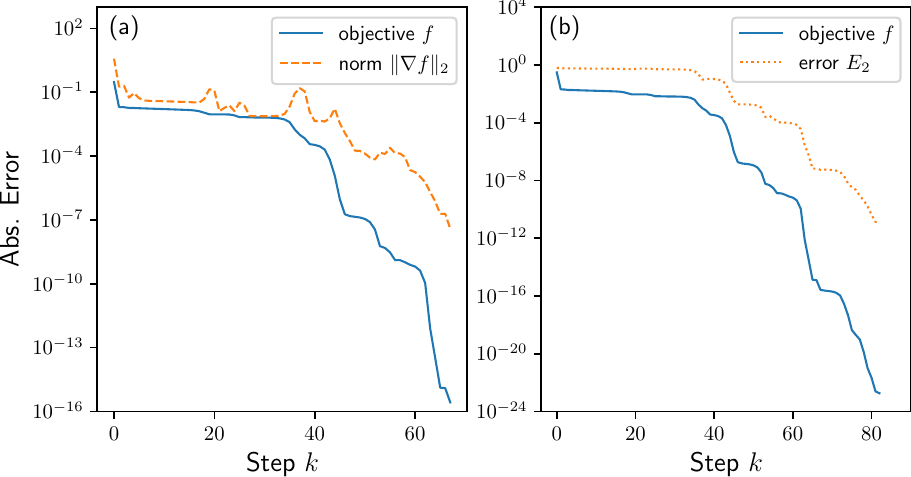}
\end{center}
\caption{\label{fig:conv_embryo} Convergence of SR1 trust-region method in case of embryo model for $s=-10$: (a) objective function and gradient norm for initial position from (\ref{eq:def_bY_embryo}) and (b) objective function and maximum Euclidean error from (\ref{eq:def_bY_embryo}).} 
\end{figure}

First, we repeat the same procedure with $s=-10$ as for the symmetric case, and Figure \ref{fig:snaps_foam} shows snapshots of the surface $S_2$ at a few steps $k$. In this case, it seems that both variables $b$ and $Y$ are updated from the beginning and then the rotational part is further optimized subsequently. Whereas the case of icosahedra does not converges for (\ref{eq:def_bY_embryo}) with $s=-10$, the suggested algorithm applied to the embryo model succeeds in finding the inverse $\varphi(X)^{-1}$, which is also shown in Fig.\ref{fig:conv_embryo}(a) as substantial decays in the objective function and the gradient norm after two flat intervals. For reference, this optimization takes 67 steps (88.9 sec) in total. \par

Now we move on to discussion on numerical results from a more geometric perspective. That is, even though we have shown that the suggested algorithm is capable of finding zeros of the objective function (\ref{eq:def_objective}) with reasonable initial guesses, it is still nontrivial to estimate to what extent such decay of the weak distance is reflected on more geometric concepts. To answer this question, we evaluate the maximum $l^2$ norm between nodes constituting the two surfaces and plot them as a function of the step $k$. More formally, suppose that the triangulated surface $S_1$ consists of the nodes $\{p_{1,i}\}$ and $S_2$ with $\{p_{2,i}\}$ is created by applying an isometry to $\{p_{1,i}\}$. Then, at each step $k$, we evaluate the metric-based error
\begin{equation}
\label{eq:def_max_l2}
E_2 = \max_{i} \|p_{1,i} - p_{2,i}\|_2,
\end{equation}
and plot it (dotted line) in Fig.\ref{fig:conv_embryo}(b) along with the corresponding value of $f$. Here, the optimization continues until it achieves $\|\nabla f\|_2< 10^{-11}$ and takes 82 steps (118 sec). It is clearly shown that behavior of the error $E_2$ is closely linked to that of the objective function, and it eventually reaches at the level of $10^{-11}$ as the function $f$ decreases to $f\approx10^{-23}$. Hence, it is reasonable to conclude that the weak distance (\ref{eq:def_distance_Sobolev}) actually contains geometric information even for $s=-10$, which partly justifies the choice of the parameter $s$ in the previous work \cite{Koga2024}. However, it also should be noted that such a tight link between the two quantities may be peculiar to optimization on the group $SE(3)$ and further investigations in a more flexible setting is desirable. Towards the goal, the surface reconstruction problem based on a Fourier-type loss function \cite{IcSiDeTu2025} is a promising candidate, and we leave it for future works.

\section{Conclusions}\label{conclusion}
In this paper, we have considered local optimization of the weak distance between two compact surfaces in $\mathbb{R}^3$ as the associated surface measures on the special Euclidean group $SE(3)$. By applying an isometry to one surface, the distance can be regarded as a function on the space of isometries, and for appropriate exponents of the inhomogeneous Sobolev norm (\ref{eq:sobolev_2}), the second power of the function (\ref{eq:def_objective}) acquires sufficient differentiability that allows to search for its local minima in a derivative-based framework. In particular, essential part of computing the gradient of the objective function is reduced to evaluations of three additional Fourier-type transforms (\ref{eq:def_fourier_weight}), each of which can be efficiently approximated by NUFFT. In numerical experiments, we have observed convergence and accuracy of the SR1 trust-region method in a few root-finding problems for triangulated surfaces, and also have verified that the maximum value of the Euclidean distance between the nodes closely follows behavior of the objective function, which implies a tight link between the weak distance and more geometric quantities. As concluding remarks, we mention three important directions for further developments. \par

First, as shown in Section \ref{results}, the optimization problem on the group $SE(3)$ is generally nonconvex, and the suggested algorithm does not necessarily return the global minimizer if the initial guess is randomly chosen. Therefore, in applications, it is crucial to first seek initial guesses for which the local algorithm is likely to converges to a global solution, and this step is more significant in partial alignment where one of the two geometric data is fragmented. In this context, while exhaustive search using FFT in VESPER \cite{HaTeChChKi2021} is useful in practice, the global algorithm built upon the unbalanced Gromov-Wasserstein theory \cite{RiZhCoChDu2025} succeeds in extracting rough approximations of isometries as solutions of partial alignment problems for EM density maps and may be combined with the local algorithm in the present work. This idea is somewhat similar to one in \cite{RiZhCoChDu2025} where outputs from the global search are refined by the \texttt{fitmap} command in ChimeraX \cite{PeGoHuMeCoCrMoFe2021}. However, \texttt{fitmap} performs the steepest gradient ascent for maximizing correlations between data, and it should be advantageous to replace the first-order method with the suggested algorithm that employs the trust-region method and hence is more reliable.\par

Second, we have pursued accuracy to reveal potentials and limitations of the distance (\ref{eq:def_distance_Sobolev}) and its numerical optimization, and further pushing down computational time of the algorithm is yet desirable in practical applications. In addition to acceleration at the hardware level such as parallelization on multi-GPU environments and switching to the single-precision floating-point arithmetics (i.e., FP32), promising options at the early stages of optimization include coarsening numerical quadratures on both spatial and frequency domains, increasing the tolerance in NUFFT, and adaptively controlling the exponent $s$ of the Sobolev norm. In particular, effectiveness of the third strategy is partly shown in Section \ref{results}, where increasing the regularity $s$ gives improvements to convergence of the algorithm. \par

Last, the present work only assumes that given data can be regarded as Borel measures of compact support and their Fourier transforms are computationally accessible. Therefore, we should be able to extend the main ideas to geometric data of other dimensions, and apply it to, for example, aligning EM density maps as compactly supported nonnegative functions and mapping protein conformations as weighted point clouds into the quotient space of a high-dimensional Euclidean space by the orthogonal group $O(3)$ \cite{DiEsLeOkSch2024}. We hope to address these problems in future works.


\section*{Acknowledgement}
The author would like to thank Lee Ching-pei for suggesting the use of the trust-region method in this project, and Alex Barnett for valuable comments on error analysis of the trapezoidal rule.

\bibliographystyle{siamplain}
\bibliography{ref}

\end{document}

%% file: shared.tex
\usepackage{lipsum}
\usepackage{amsfonts}
\usepackage{graphicx}
\usepackage{epstopdf}
\usepackage{algorithmic}
\ifpdf
  \DeclareGraphicsExtensions{.eps,.pdf,.png,.jpg}
\else
  \DeclareGraphicsExtensions{.eps}
\fi

\newsiamremark{remark}{Remark}
\newsiamremark{hypothesis}{Hypothesis}
\crefname{hypothesis}{Hypothesis}{Hypotheses}
\newsiamthm{claim}{Claim}
\newsiamremark{fact}{Fact}
\crefname{fact}{Fact}{Facts}

\headers{Local optimization of Weak Distance between compact surfaces}{K. Koga}

\title{Local optimization of weak distance between compact surfaces on special Euclidean group\thanks{Submitted June 5, 2026.
\funding{This work was partially supported by, JST ACT-X JPMJAX2106, JST PREST JPMJPR25K2, and JST CREST JPMJCR22P5.}}}

\author{Kazuki Koga\thanks{Department of Mathematical and Computing Science, Institute of Science Tokyo, Meguro-ku, Tokyo, Japan.}}

\usepackage{amsopn}
